\documentclass[10pt]{article}
\usepackage[cp1251]{inputenc}
\usepackage[english]{babel}
\usepackage{amssymb}
\usepackage{amsmath}
\usepackage{xcolor}
\usepackage{fancyhdr} 
\usepackage{graphicx}

\begin{document}
\begin{center}
\bf {\Large Automata as $p$-Adic Dynamical Systems}
   \end{center}   

\begin {center}{Livat B. Tyapaev}
\end{center}
\begin{center}{e-mail: TiapaevLB@info.sgu.ru}
\end{center}
\begin{center}{National Research Saratov State University}
\end{center}
\begin{abstract}
{\small
The automaton transformation of infinite words over alphabet $\mathbb F_p=\{0,1,\ldots,p-1\}$, where $p$ is a prime number, coincide with the continuous transformation (with respect to the $p$-adic metric) of a ring $\mathbb Z_p$ of $p$-adic integers.
The objects of the study are non-Archimedean dynamical systems generated by automata mappings on the space $\mathbb Z_p$. Measure-preservation (with the respect to the Haar measure) and ergodicity of such dynamical systems plays an important role in cryptography (e.g. for pseudo-random generators and stream cyphers design). The possibility to use $p$-adic methods and geometrical images of automata allows to characterize of a transitive (or, ergodic) automata. We investigate a measure-preserving and ergodic mappings associated with synchronous and asynchronous automata. We have got criterion of measure-preservation for an $n$-unit delay mappings associated with asynchronous automata. Moreover, we have got a sufficient condition of ergodicity of such mappings.}
\end{abstract}

\section{Introduction}
Speaking about a (synchronous) automaton $\mathfrak A$, we always understand the letter-to-letter transducer (not necessarily with a finite number of states) with fixed initial state. Input and output alphabets $\mathbb F_p=\{0,1,\ldots,p-1\}$ of an automaton consist of $p$ symbols, where $p$ is a prime number. The automaton naturally defines the mapping of the set $\mathbb Z_p$ of all (one-sided) infinite sequences in $\mathbb Z_p$. As is known, such maps are called {\it deterministic functions (or, automata functions)}. The set $\mathbb Z_p$ is naturally endowed with the structure of the ring of a $p$-adic integers, i.e. turn into the metric space by specifying the metric $|a-b|_p$, $a,b\in\mathbb Z_p$, where $|\cdot|_p$ is the $p$-adic (that is, non-Archimedean) norm. {\it All deterministic functions $f\colon\mathbb Z_p\to\mathbb Z_p$ satisfy the $p$-adic Lipschitz condition with constant 1, i.e. $|f(a)-f(b)|_p\le|a-b|_p$ for all $a,b\in\mathbb Z_p$}, see \cite {AnKh}. In particular, all deterministic functions are continuous with respect to the $p$-adic metric. For example, all functions defined by polynomials with $p$-adic integers coefficients (in particular, rational integers) are deterministic. Note, that in the case of $p=2$, the standard processor commands, such as AND, OR, XOR, NOT, etc. naturally extend to functions from $\mathbb Z_2\times \mathbb Z_2$ in $\mathbb Z_2$. Hence, in particular, it follows that all functions obtained with the compositions of arithmetic and coordinate-wise logical operations of the processor can be considered as $2$-adic deterministic functions.
By the above, $p$-adic analysis can prove to be a very effective ``analytical"  tool for research the properties of deterministic functions and the behavior of automata.

\section{Transitivity of automata}
We consider only automata, where every state $s$ from the set $\mathcal S$ of all internal states of an automaton is accessible from the initial state $s_0$, that is, the automaton will pass from $s_0$ to $s$ by reading of some finite input word. Every automaton $\mathfrak A$ defines a family of automata $\mathfrak A_s$, where $\mathfrak A_s$ differs from $\mathfrak A$ only in that it has a different initial state, $s$ instead of $s_0$. An automaton $\mathfrak A$ is said to be {\it transitive (or, ergodic)} if the family of deterministic functions defined by the family of automata $\mathfrak A_s$, $s\in\mathcal S$ is transitive on every set $\mathbb Z/p^n\mathbb Z$ for all $n=1,2,\ldots$ (that is, for all finite words $u$ and $v$ of the same length there exists an automaton $\mathfrak A_s$ that transforms $u$ into $v$). We associate the automaton $\mathfrak A$ with the closure  $\mathcal E_{\mathfrak A}$ of all points of the form $(\frac{u}{p^n},\frac{v}{p^n})\in[0,1]\times[0,1]$ in the topology of the Euclidean plane, where $u$ is the input word of length $n$, and $v$ is the corresponding output word, $n=1,2,\ldots$. The set $\mathcal E_{\mathfrak A}$ is measurable with respect to the Lebesgue measure $\lambda$. The following ``law 0-1" holds: {\it for any automaton $\mathfrak A$, $\lambda(\mathcal E_{\mathfrak A})=0 $ or $\lambda(\mathcal E_{\mathfrak A})=1$. Moreover, $\lambda (\mathcal E_{\mathfrak A})=1$ if and only if $\mathfrak A$ is transitive}, see \cite {AnKh}.

Let us enumerate symbols $\alpha_i$ of the alphabet $\mathbb F_p$ by natural numbers $c(\alpha_i)\in\{1,2,\ldots p\}$, and we associate with the word $u=\alpha_{k-1}\ldots\alpha_1\alpha_0$ over the alphabet $\mathbb F_p$ the rational number $\overrightarrow u=c(\alpha_0)+c(\alpha_1)\cdot(p+1)^{-1}+\ldots+c(\alpha_{k-1})\cdot(p+1)^{-(k-1)}$. We consider all points of the Euclidean square $[1,p+1]\times[1,p+1]\subset \mathbb R^2$ of the form $(\overrightarrow u, f(\overrightarrow u))$, where $u$ runs through all finite words over $\mathbb F_p$. The set of these points $\Omega(\mathfrak A)$ is called {\it a geometrical image (or, graph)} of an automaton $\mathfrak A$, see \cite {T4,T5}. For every state $s\in \mathcal S$ of an automaton $\mathfrak A$ we associate a map $R_s\colon \mathbb F_p\to\mathbb F_p$ that transforms input symbol into output symbol. If we consider an automaton $\mathfrak A$ and the family $\{R_s\colon s\in\mathcal S\}$, then a correspondence for every state $s\in\mathcal S$ of a some map $R_s$ creates a new automaton $\mathcal B$ (and a set $K(\mathfrak A)$ of an automata that is constructed this way). {\it Then automaton $\mathfrak A$ is transitive if and only if there exists the automaton $\mathfrak B\in K(\mathfrak A)$ and such a geometrical images $\Omega(\mathfrak A)$, $\Omega(\mathfrak B)$, that are affine equivalents} \cite {T10, T12}.

\section{Dynamical systems}
A \textit {dynamical system} on a measurable space $\mathbb S$ is understood as a triple $(\mathbb S,\mu,f)$, where $\mathbb S$ is a set endowed with a measure $\mu$, and $f\colon\mathbb S\to \mathbb S$ is a measurable function. A dynamical system is also topological since configuration space $\mathbb S$ are not only measurable space but also metric space, and corresponding transformation $f$ are not only measurable but also continuous. 

We consider a $p$-adic dynamical system $(\mathbb Z_p, \mu_p, f)$ \cite{A2}. The space $\mathbb Z_p$ is equipped with a natural probability measure, namely, the Haar measure $\mu_p$ normalized so that $\mu_p(\mathbb Z_p)=1$. Balls $B_{{p}^{-k}}(a)$ of nonzero radii constitute the base of the corresponding $\sigma$-algebra of measurable subsets, $\mu_p(B_{{p}^{-k}}(a))=p^{-k}$. A measurable mapping $f\colon\mathbb Z_p\to \mathbb Z_p$ is called {\it measure-preserving} if $\mu_p(f^{-1}(S))=\mu_p(S)$ for each measurable subset $S\subset\mathbb Z_p$.
A measure-preserving map $f$ is said to be \textit{ergodic} if for
each measurable subset $S$ such that $f^{-1}(S)=S$ holds either $\mu_p(S)=1$ or $\mu_p(S)=0$; so ergodicity of the map $f$ just means that $f$ has no proper invariant subsets; that is, invariant subsets whose measure is neither $0$ nor $1$.

We can consider an automaton function $f\colon \mathbb Z_p\to\mathbb Z_p$ as an algebraic dynamical system on a measurable and a metric space $\mathbb Z_p$ of $p$-adic integers, which, actually, is a profinite algebra with the structure of an inverse limit: The ring $\mathbb Z_p$ is an inverse limit of residue rings $\mathbb Z/p^k\mathbb Z$, $k=1,2,3\ldots$. As any profinite algebra can be endowed with a metric and a measure, it is reasonable to ask what continuous with respect to the metric transformations are measure-preserving or ergodic with respect to the mentioned measure. Besides, the same question can be asked in the case of mappings for asynchronous automata.

\section{Measure-preserving and ergodic an $n$-unit delay mappings}
We assume that an {\it asynchronous automaton (letter-to-word transducer)} $\mathcal C$ works in a framework of discrete time steps. The transducer reads one symbol at a time, changing its internal state and outputting a finite sequence of symbols at each step. {\it Asynchronous transducers are a natural generalizations of synchronous transducers}, which are required to output exactly one symbol for every symbol read. A mapping $f_{\mathcal C}\colon\mathbb Z_p \to \mathbb Z_p$ is called {\it $n$-unit delay} whenever given an asynchronous automaton $\mathcal C$ translated infinite input string $\alpha=\ldots \alpha_n\alpha_{n-1}\ldots \alpha_1\alpha_0$ (viewed as $p$-adic integer) into infinite output string $\beta=\ldots \beta_{n+1}\beta_{n}$ (viewed as $p$-adic integer). An $n$-unit delay transducer produces the some output $n$ times unit later. Note that usually the term $n$-unit delay is used in a narrower meaning, cf. \cite {Gregorchuk} when  $n$-unit delay transducer is defined by finite automaton, that produces no output for the first $n$ times slots; after that, the automaton outputs the incoming words without changes. {\it An $n$-unit delay mapping $f_{\mathcal C}\colon\mathbb Z_p\to \mathbb Z_p$ is continuous on $\mathbb Z_p$} \cite {T12}.

Let $F_k$ be a reduction of function $f$ modulo $p^{n\cdot(k-1)}$ on the elements of the ring $\mathbb Z/p^{n\cdot k}\mathbb Z$ for $k=2,3,\ldots$. The following criterion of measure-preservation for $n$-unit delay mappings is valid: {\it An $n$-unit delay mapping $f\colon\mathbb Z_p \to \mathbb Z_p$ preserves the measure if and only if the number $\#F_k^{-1}(x)$ of $F_k$-pre-images of the point $x\in \mathbb Z/p^{n\cdot(k-1)}\mathbb Z$ is equal $p^n$, $k=2,3,\ldots$} \cite {T11,T12}.

A point $x_0\in \mathbb Z_p$ is said to be a \textit{periodic point} if there exists $r\in\mathbb N$ such that $f^r(x_0)=x_0$. The least $r$ with this property is called the \textit{length} of period of $x_0$. If $x_0$ has period $r$, it is called an \textit{$r$-periodic point}. The orbit of an $r$-periodic point $x_0$ is $\{x_0, x_1,\ldots,x_{r-1}\}$, where $x_j=f^j(x_0)$, $0\le j\le r-1$. This orbit is called an \textit{$r$-cycle}. Let $\gamma(k)$ be an $r(k)$-cycle $\{x_0,x_1,\ldots,x_{r(k)-1}\}$, where $$x_j=(f\mod p^{k\cdot n})^j(x_0),$$ $0\le j\le r(k)-1$, $k=1,2,3,\ldots$. The following condition of ergodicity holds: {\it A measure-preserving an $n$-unit delay mapping $f\colon\mathbb Z_p\to\mathbb Z_p$ is ergodic if $\gamma(k)$ is an unique cycle, for all $k\in \mathbb N$} \cite {T13}.

By Mahler's Theorem, any continuous function $f\colon \mathbb Z_p\to \mathbb Z_p$ can be expressed in the form of a uniformly convergent series, called its \emph {Mahler Expansion}:
$$f(x)=\sum_{i=0}^{\infty}a_i\binom{x}{i},$$
where $a_i\in\mathbb  Z_p$, $i = 0,1,2,\ldots$, and
$$\binom{x}{i}=\frac{x(x-1)\cdots(x-i+1)}{i!}$$
for $i = 1,2,\ldots$; 
$$\binom{x}{0}=1,$$ by the definition.
For an $n$-unit delay mapping, $n\in \mathbb N$, we gets next theorem.

Theorem 1. {\it A function $f\colon \mathbb Z_p\to \mathbb Z_p$ represented by Mahler expansion
$$f(x)=\sum_{m=0}^{\infty}a_m\binom{x}{m},$$ where $a_m\in\mathbb Z_p$, $m=0,1,2\ldots$,
is an $n$-unit delay if and only if
$$|a_i|_p\le p^{-\lfloor \log_{p^n}i\rfloor+1}$$
for all $i\ge 1$}.

The following criterion of measure-preservation for $n$-unit delay mappings is valid.

Theorem 2. {\it An $n$-unit delay mapping $f\colon\mathbb Z_p\to\mathbb Z_p$ is measure-preserving whenever the following conditions hold  simultaneously:
\begin{enumerate}
\item $a_i\not \equiv 0\pmod p$ for $i=p^n$;
\item $a_i\equiv 0 \pmod {p^{\lfloor \log_{p^n}i\rfloor}}$, $i>p^n$. 
\end{enumerate}}

The following condition of ergodicity holds.

Theorem 3. {\it An $n$-unit delay mapping $f\colon\mathbb Z_p\to\mathbb Z_p$ is ergodic on $\mathbb Z_p$ whenever the following conditions hold  simultaneously:
\begin{enumerate}
\item $a_1+a_2+\ldots+a_{p^n-1}\equiv 0 \pmod p$;
\item $a_i\equiv 1\pmod p$ for $i=p^n$;
\item $a_i\equiv 0 \pmod {p^{\lfloor \log_{p^n}i\rfloor}}$, $i>p^n$. 
\end{enumerate}}

Given an automaton $\mathcal C$, consider the corresponding an $n$-unit delay mapping $f=f_{\mathcal C}\colon\mathbb Z_p\to\mathbb Z_p$. Let $E_k(f)$ be a set
of all the following points $e_k^f (x)$ of Euclidean unit square $\mathbb I^2=[0, 1]\times[0, 1]\subset
\mathbb R^2$ for $k=1,2,3,\ldots$:
$$e_k^f (x)=\Bigl(\frac{x\bmod p^{n+k}}{p^{n+k}},\frac{f(x)\bmod p^{k}}{p^{k}}\Bigr),$$
where $x\in\mathbb Z_p$.
Note that $x \bmod p^{n+k}$ corresponds to the prefix of length $n+k$ of the infinite word $x\in\mathbb Z_p$, i.e., to the input word of length $n+k$ of the automaton $\mathcal C$; while $f(x) \bmod p^k$ corresponds to the respective output word of length $k$. Denote $\mathcal E(f)$ the closure of the set $E(f)=\bigcup_{k=1}^{\infty}E_k(f)$ in the topology of real plane $\mathbb R^2$. As $\mathcal E(f)$ is closed, it is measurable with respect to the Lebesgue measure on real plane $\mathbb R^2$. Let $\lambda(f)$ be the Lebesgue measure of $\mathcal E(f)$.

Theorem 4. {\it For an $n$-unit delay mapping $f\colon\mathbb Z_p\to\mathbb Z_p$ the closure $\mathcal E(f)$ is nowhere dense in $\mathbb I^2=[0,1]\times[0,1]\subset R^2$ (equivalently, $\lambda(f)=0$)}.


\begin{thebibliography}{30}

\bibitem{AnKh} V. Anashin and A.\;Khrennikov, {\it Applied Algebraic Dynamics}, de Gruyter Expositions
in Mathematics (de Gruyter GmbH~\&~Co., Berlin--N.Y., 2009).

\bibitem {A2} V. Anashin, “Ergodic transformations in the space of $p$-adic integers,” Proc. Int. Conf. on p-adic Mathematical Physics, AIP Conference Proceedings 826, 3–24 (2006).

\bibitem{Gregorchuk} R.\,I.\;Grigorchuk,\;V.\,V.\;Nekrashevich, and V.\,I.\;Sushchanskii, ``Automata, dynamical systems, and groups," Proc. Steklov Math. Inst. {\bf 231}, 128–203 (2000).

\bibitem {T4} L.\,B.\;Tyapaev, ``The geometrical model of behavior of automata and their indistinguishability," Mathematics, Mechanics {\bf 1}, 139-143 (1999) [in Russian].

\bibitem {T5} L.\,B.\;Tyapaev, ``Solving some problems of automata behaviour analysis,'' Izv. Saratov Univ. (N.S.), Ser. Math. Mech. Inform. {\bf 6}:1-2, 121-133 (2006).

\bibitem {T10} L.\,B.\;Tyapaev, ``Transitive families of automata mappings," in {\it Proceedings of the 9th International Conference on Discrete Models in the Theory of Control Systems} (May 20-22 2015, Moscow), eds. V.\,B.\;Alekseev,\;D.\,S.\;Romanov,\;B.\,R.\;Danilov, 244-247 (Lomonosov Moscow State University, Maks Press, 2015) [in Russian].

\bibitem {T11} L.B. Tyapaev, ``Measure-preserving and ergodic asynchronous automata mappings," in {\it Proceedings of the 12th International Workshop on Discrete Mathematics and it Applications} (June 20-26 2016, Moscow), edited by O.M. Kasim-Zade, 398-400 (published by Faculty of Mechanics and Mathematics, Lomonosov Moscow State University, 2016) [in Russian].

\bibitem {T12} L.B. Tyapaev, ``Transitive families and measure-preserving an $n$-unit delay mappings," in {\it Proceedings of the International Conference on Computer Science and Information Technologies} (June 30-July 2 2016, Saratov), 425-429 (Publishing Center Nauka, Saratov, 2016).

\bibitem {T13} L.B. Tyapaev, ``Ergodic automata mappings with delay," in {\it Proceedings of the International Conference on Problems of theoretical cybernetics} (June 19-23 2017, Penza), edited by Yu. I. Zhuravlev, 242-244, (Moscow, Maks Press, 2017) [in Russian].

\end{thebibliography}
\end{document}